\documentclass[11pt, a4paper]{amsart}
 
\usepackage[euler-digits]{eulervm}

\usepackage{amsfonts, amsthm, amssymb, amsmath, stmaryrd}
\usepackage{mathrsfs,array}
\usepackage{eucal,fullpage,times,color,enumerate,accents}
\usepackage{url}

\usepackage{color}
\usepackage{mathrsfs}
\usepackage{amssymb}
\usepackage{tikz}
\usepackage{tikz-cd}
\usepackage{bm}
\usepackage{enumerate}
\usetikzlibrary{calc}

\DeclareMathAlphabet{\mathpzc}{OT1}{pzc}{m}{it}
\makeatletter
\def\blfootnote{\xdef\@thefnmark{}\@footnotetext}
\makeatother

\title{\raisebox{0.1in}{Superabundant curves and the Artin fan}}
\author{Dhruv Ranganathan}
\date{\today}

\usepackage[backref]{hyperref}
\hypersetup{
  colorlinks   = true,          %Colors links instead of ugly boxes
  urlcolor     = purple,          %Color for external hyperlinks
  linkcolor    = blue,          %Color of internal links
  citecolor   = purple             %Color of citations
}

\address{Department of Mathematics, Yale University}
\email{dhruv.ranganathan@yale.edu}

\newtheorem*{mainthm}{Main Theorem}
\newtheorem{theorem}{Theorem}

\newtheorem*{question}{Question}

\newtheorem{example}[theorem]{Example}
\newtheorem{quasi-theorem}[theorem]{Quasi-Theorem}
\newtheorem{blank remark}[theorem]{}

\newtheorem{rem1}[theorem]{Remark}
\newenvironment{remark}{\begin{rem1}\em}{\end{rem1}}

\newtheorem{not1}[theorem]{Notation}

\newcommand{\CC} {{\mathbb C}}          
            
\newcommand{\NN} {{\mathbb N}}		
\newcommand{\PP}{\mathbb{P}}         
\newcommand{\QQ} {{\mathbb Q}}		
\newcommand{\RR} {{\mathbb R}}		
\newcommand{\ZZ} {{\mathbb Z}}

\newcommand{\Hom}{\operatorname{Hom}}

\DeclareMathOperator{\val}{val}

\DeclareMathOperator{\spec}{Spec}

\def\cA{{\mathpzc A}}

\def\fM{\mathfrak{M}}

\def\fm{\mathfrak{m}}

\def\trop{\mathrm{trop}}
\def\an{\mathrm{an}}

\begin{document}
\vspace{-0.25in}
\pagestyle{plain}
\maketitle

\begin{abstract}
We prove that every balanced $1$-dimensional polyhedral complex arises as the tropicalization of a smooth curve over a non-Archimedean field mapping to a toric Artin fan, namely the quotient of a toric variety by its dense torus. 
\end{abstract}

\blfootnote{This research was partially supported by NSF grant CAREER DMS-1149054 (PI: Sam Payne).}

\section{Introduction}
Let $K$ be a non-Archimedean field and $T$ an algebraic torus with character lattice $M$. The Berkovich analytic space $T^{\an}$ is a space of real valuations on $K[M]$. There is a natural continuous map $\trop:T^{\an}\to \Hom(M,\RR)$, obtained by restricting valuations to the character lattice. When $C$ is a curve in $T$, the image of $C^{\an}$ under $\trop$ is a metric graph, embedded in $\Hom(M,\RR)$ as a rational polyhedral complex $\mathscr P$ of dimension $1$. Each edge $e$ carries a weight, equal to the number of components, counted with multiplicity, of the initial degeneration defined by any point $p$ on the interior of that edge~\cite{MS14}. The weights satisfy the so-called balancing condition at vertices -- the sum of the outgoing slopes at every vertex is zero. Define an \textit{embedded tropical curve} to be a balanced, connected, rational polyhedral complex of dimension $1$ in a vector space. Consider the inverse problem.

\begin{question}
Given an embedded tropical curve $\mathscr P$ in $\Hom(M,\RR)$, does there exist an algebraic curve $C$ in $T$ whose tropicalization is $\mathscr P$?
\end{question}

\begin{figure}
\includegraphics[scale=0.2]{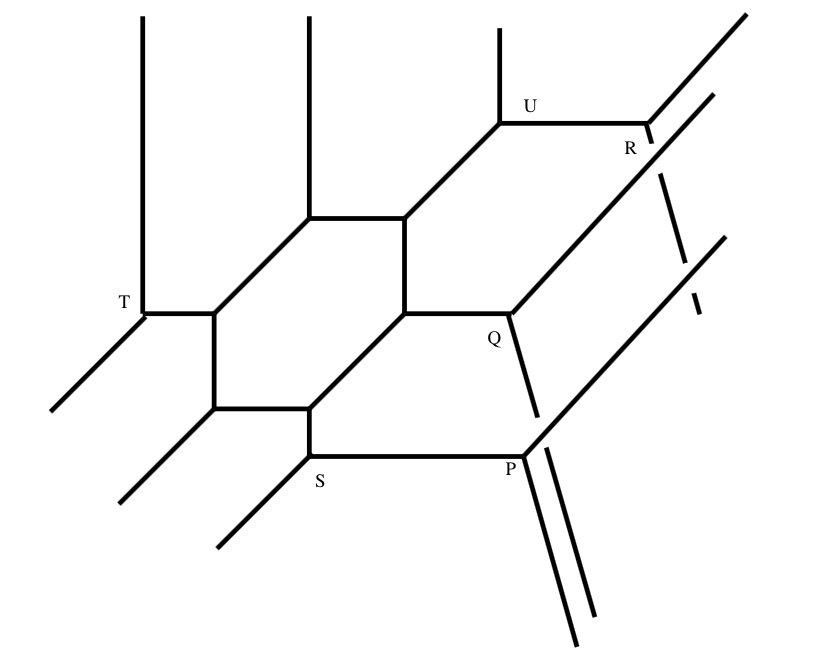}
\caption{The Mikhalkin--Speyer example of a genus $1$ tropical curve that cannot be lifted to an algebraic curve~\cite[Figure 5.1]{Sp-thesis}. This graph fails to fulfill Speyer's well-spacedness condition~\cite[Theorem 3.3]{Sp07}.}
\label{speyercurve}
\end{figure}

\noindent
The answer to this question is positive when $M$ has rank $2$, or when the genus of $\mathscr P$ is zero, see~\cite{Mi03,NS06}. However, there exist higher genus tropical curves in $\RR^n$ that are not the tropicalizations of algebraic curves in $(K^*)^n$. For instance Mikahlkin and Speyer~\cite{Sp-thesis, Sp07} constructed an example of a genus $1$ tropical curve in $\RR^3$ that does not arise as a tropicalization, depicted in Figure~\ref{speyercurve}. This is due to the phenomenon of \textit{superabundance}. That is, of tropical curves in $\RR^n$ having spaces of deformations that are strictly larger than the expected dimension. See~\cite[Definition 2.2]{Mi03} and the discussion in~\cite[Section 1]{Kat12}. A number of authors have studied this lifting problem in a variety of contexts, see~\cite{CDMY,CFPU,Kat12, Ni15,NS06}.

We consider the question of whether these superabundant tropical curves are tropicalizations in some more general sense. Recently, Ulirsch~\cite{U14b} has shown that the tropicalization map for a proper toric variety $X(\Delta)$ coincides with the map from $X(\Delta)^{\an}$ to the analytic stack quotient $[X(\Delta)^{\an}/T^{\an}_\circ]$, where, $T^{\an}_\circ$ is the non-Archimedean analytic compact torus, consisting of valuations that are identically $0$ on the character lattice of $T$. The topological space underlying the stack $[X(\Delta)^{\an}/T^{\an}_\circ]$ is canonically identified with the extended tropicalization of $X(\Delta)$ defined by Kajiwara~\cite{Kaj08} and Payne~\cite{Pay09}.

Let $\mathscr P$ be an embedded tropical curve. Denote by $\overline{\mathscr P}$ be the compactification of $\mathscr P$ obtained by adding a point at infinity to compactify each unbounded edge.

\begin{mainthm}\label{mainthm}
There exists a smooth curve over a non-Archimedean field $K$ and a proper toric variety $X(\Delta)$, such that the extended polyhedral complex $\overline{\mathscr P}$ coincides with the topological image of a map of analytic stacks
\[
C^{\an}\to [X(\Delta)^{\an}/T^{\an}_\circ].
\]
The weights on edges of $\mathscr P$ are recovered via the expansion factors of the map above. Moreover, $\Delta$ can be chosen to be any complete fan whose one-skeleton supports the recession fan of $\mathscr P$. 
\end{mainthm}

In other words, the realizability problem for tropical curves amounts precisely to lifting this map $C^{\an}\to[X(\Delta)^{\an}/T^{\an}_\circ]$ to a map $C^{\an}\to X(\Delta)^{\an}$. 

\begin{remark}
Geometrically, maps from a curve $C$ to an Artin fan can be understood as collections of line bundles and sections on $C$. For instance, when $X(\Delta) = \PP^n$, a map to $[\PP^n/T]$ is a collection of $n+1$ line bundles with sections $(L_0,s_0),\ldots, (L_n,s_n)$  on $C$, such that, at every point of $C$, at least one section $s_j$ is nonzero. When all the line bundles coincide, one obtains a map to $\PP^n$. 
\end{remark}

The approach to the proof of the main theorem in similar in spirit to the techniques developed by Nishinou--Siebert~\cite{NS06}, and extended by Cheung--Fantini--Park--Ulirsch~\cite{CFPU} and Nishinou~\cite{Ni15}. However, while toric degenerations of toric varieties play a central role in the cited works, we work entirely within the framework of logarithmic stable maps, as developed in the series of recent papers~\cite{AC11,Che10,GS13}. As a result, we avoid the ``expansion of target'' via toric degenerations, and instead rely heavily on the log structure of the stack $[X(\Delta)/T]$. The existence of algebraic moduli spaces for logarithmic prestable maps also allows us to handle tropical curves with edge lengths that are not necessarily rational. 

We freely use the Fontaine--Illusie--Kato theory logarithmic geometry in this note, and refer the reader to the surveys~\cite{ACGHOSS,ACMUW} and K. Kato's article~\cite{Kat89} for an introduction.

\subsection*{Terminology} The \textit{star} around a polyhedron $P$ of a polyhedral complex $\mathscr P$ is the fan of all polyhedra in $\mathscr P$ containing $P$, as defined in~\cite[Definition 2.3.6]{MS14}. The recession fan of a polyhedral complex $\Sigma$ is the collection of cones of unbounded directions, as defined in~\cite[Section 3.5, p.133]{MS14}.

\section{Analytification, tropicalization, and artin fans}

\subsection{Analytification} Let $K$ be a field complete with respect to a rank-$1$ valuation $\nu$. Throughout we will assume that $K$ is an extension of $\CC$ and that $\nu$ induces the trivial valuation on $\CC$. The valuation ring and maximal ideal will be denoted $R$ and $\fm$ respectively.

Let $X$ be a finite type $K$-scheme. The \textit{Berkovich analytification} $X^{\an}$ is constructed as a locally ringed space in~\cite{Ber90}. For affine $X = \spec(A)$, one considers the following set of ring valuations on $A$,
\[
X^{\an} = \{\val:A\to \RR\sqcup \{\infty\}:\val|_K = \nu\}.
\]
The set $X^{\an}$ is given the weak topology for the evaluation functions
\begin{eqnarray*}
ev_f: X^{\an}&\to& \RR\sqcup \{\infty\},\\
\val_x&\mapsto& \val_x(f).
\end{eqnarray*}
For arbitrary $X$, one uses the above construction on affines and glues the resulting pieces, giving rise to a functor
\[
(-)^{\an}: \mathbf{Schemes}/K\to \mathbf{An.Spaces}/K.
\]
In particular, every point of $X^{\an}$ can be represented by a map $\spec(L)\to X$, where $L$ is a valued field extending $K$, compatible with the valuation on $K$. 

\subsection{Generic fiber} Assume the valuation on $K$ is nontrivial. The analytic space $X^{\an}$ is Hausdorff (resp. compact) if and only if the scheme $X$ is separated (resp. proper)~\cite[Theorem 3.4.8]{Ber90}. Since the algebraic stacks $[X(\Delta)/T]$ are never separated, to relate them to tropicalizations we need a variation of this analytification, which is known as \textit{Raynaud's generic fiber}~\cite[Section 7.4]{Bo05}. 

Given a flat $R$-scheme $\mathscr Z$ of finite type, Raynaud's generic fiber, usually denoted $\widehat {\mathscr Z}_\eta$, is a Berkovich analytic space over $K$, associated to the formal completion of $\mathscr Z$ along the maximal ideal of $R$. We work exclusively with integral $R$-schemes, as opposed to formal schemes. In order to avoid confusion we define $\mathscr Z^{\an}_\circ: = \widehat {\mathscr Z}_\eta$. 

The generic fiber can be described as follows. Let $\mathscr Z$ be an affine finite-type scheme over $\spec(R)$ with generic fiber $Z$. In this case $\mathscr Z^{\an}_\circ$ is the compact analytic domain in $Z^{\an}$ consisting of points represented by maps $\spec(L)\to Z$, that extend to maps $\spec(R_L)\to \mathscr Z$ from the valuation ring of $L$. This construction can be extended to all $R$-schemes by a gluing process, yielding a functor~\cite[Section 0.3.3]{Bar96}
\[
(-)^{\an}_\circ: \mathbf{Schemes}/R\to \mathbf{An.Spaces}/K.
\]
If $\mathscr Z$ is proper, the valuative criterion of properness ensures that $\mathscr Z^{\an}_\circ$ coincides with $Z^{\an}$. In general one has a map $\mathscr Z^{\an}_\circ\to Z^{\an}$, which is injective when $\mathscr Z$ is separated over $R$. In~\cite[Section 6]{Yu14b} Yu extends the generic fiber to a functor
\[
(-)^{\an}_\circ: \mathbf{Alg.Stacks}/R\to \mathbf{An.Stacks}/K.
\]
See~\cite{U14b,Yu14b} for details on non-Archimedean analytic stacks. Note that for generic fibers and analytifications of algebraic stacks, the two notions of analytic stacks lead to the same underlying topological spaces. As a result, either theory suffices for the purposes of this article.

\begin{example}
Suppose that $M$ is a lattice and $\mathscr T = \spec(R[M])$ is a split torus over $R$, with generic fiber $T$. The space $\mathscr T^{\an}_\circ$ consists of those valuations $\val_t\in T^{\an}$ such that $\val_t(\chi^u) = 0$ for all $u\in M$. That is, $\mathscr T^{\an}_\circ$ is the non-Archimedean analytic analogue of the real torus. 
\end{example}

\noindent
As the torus $T$ over $K$ has a canonical model over $R$, we denote its Raynaud fiber by $T^{\an}_\circ$.

\subsection{$\beth$--space} When $K$ carries the trivial valuation, there is an analogue of the generic fiber construction defined by Thuillier~\cite{Thu07}:
\[
(-)^{\beth}: \mathbf{Schemes}/K\to \mathbf{An.Spaces}/K.
\]
Intuitively, this may be thought of as a generic fiber in the case where the field $K$ and its valuation ring $R$ coincide. If $X = \spec(A)$ is an affine $K$-scheme, the space $X^\beth$ is a compact analytic domain in $X^{\an}$ consisting of those multiplicative valuations $\val_x:A\to \RR\sqcup\{\infty\}$ that are nonnegative. As before, one may glue the $\beth$-spaces of affine patches to obtain a $K$-analytic space. If $X$ is separated, $X^\beth$ is a subspace of $X^{\an}$, and if $X$ is proper, $X^{\an} = X^\beth$. We refer the reader to~\cite{Thu07}, for details. In~\cite{U-thesis}, the $\beth$-space construction is extended to a functor from algebraic stacks over $K$ to analytic stacks over $K$, in analogous fashion to Raynaud's generic fiber.

\subsection{Tropicalization and the Artin fan} Let $T = \spec(K[M])$ be a torus with character lattice $M$ and dual lattice $N$. Given a point $\val_t\in T^{\an}$, one may restrict the valuation to the character lattice $M$ of $T$ to obtain a point $\trop(t)$ of $\Hom(M,\RR)$. This yields a continuous tropicalization map
\[
\trop: T^{\an}\to N_\RR.
\]
In~\cite{Pay09}, this construction is extended, replacing the torus $T$ by an arbitrary toric variety $X(\Delta)$. This yields a continuous map
\[
X(\Delta)^{\an}\to N(\Delta),
\]
where $N(\Delta)$ is a partial compactification of the vector space $N_\RR$. The tropicalization of a subvariety $Y$ of $X(\Delta)$ is defined to be the image of $Y^{\an}\hookrightarrow X(\Delta)^{\an}$ under this map.

To simplify the discussion, we henceforth assume that $\Delta$ is a complete fan. Given a toric variety $X(\Delta)$ with dense torus $T$, the quotient stack $\cA(\Delta) := [X(\Delta)/T]$ is referred to as the \textit{Artin fan} of $X(\Delta)$. The Artin fan originates from ideas in Olsson's work on the moduli space of logarithmic structures in~\cite{Ols03}. We refer the reader to~\cite{ACMUW,AW,U-thesis} for a more complete treatment of Artin fans. 

\begin{theorem}[{\cite[Theorem 1.4]{U14b}}] 
There is a natural isomorphism of extended cone complexes $\mu_\Delta: |\cA(\Delta)^{\an}_\circ|\to N(\Delta)$, making the diagram
\[
\begin{tikzcd}
\phantom{1} & \left|{\cA}(\Delta)^{\an}_\circ\right| \arrow{dd}{\mu_\Delta} \\
X(\Delta)^{\an} \arrow[swap]{dr}{\trop} \arrow{ru}{\textnormal{Stack Quotient}} & \phantom{1} \\
\phantom{1} & N(\Delta) \\
\end{tikzcd}
\]
commute. 
\end{theorem}

Here, $|-|$ is the functor associating to an analytic stack $Y$, its underlying topological space $|Y|$, as defined in~\cite[Section 5]{U14b}. 

\begin{remark}
The fact that the map from $X(\Delta)_\circ^{\an}$ to its skeleton is the quotient by the analytic group $T^{\an}_\circ$ is implicit in Berkovich's work on local contractibility~\cite{Ber99}, and in Thuillier's work in the trivially valued setting~\cite{Thu07}. Ulirsch's results allows one to enhance the topological retraction maps to analytic maps, by providing the skeleton with the structure of an analytic stack.
\end{remark}

\begin{remark}
We bring to the reader's attention an instructive analogy. If $P$ is a simple lattice polytope in $M$, the polarized complex toric variety $X(P)$ has the structure of a smooth symplectic manifold with an action of the compact torus $T_\circ = \Hom(M,\mathbb S^1)$. The quotient of the symplectic manifold $X(P)$ by $T_\circ$ coincides with the moment polytope $P$ of $X(P)$. In fact, the moment polytope $P$ and the Kajiawara--Payne extended tropicalization of $X(P)$, giving $\CC$ the trivial valuation, are isomorphic as abstract polytopes. See~\cite[Remark 3.3]{Pay09}.
\end{remark}

\section{Proof of the Main Theorem}

The main result is proved using the geometry of the moduli space of logarithmic prestable maps to the Artin fan. However, we give first an intuitive explanation of why one expects such a result to hold. When realizing tropical curves using a map to a toric variety $X(\Delta)$, one first constructs a logarithmic map from a ``would be'' special fiber of a degenerating curve to $X(\Delta)$. This curve is chosen to have dual graph equal to the underlying graph of the expected tropicalization, and the lengths, edge multiplicities, and edge directions, are encoded in the logarithmic structure. Such a map is constructed over the logarithmic base $\spec(\NN\to \CC)$, i.e. the special point in the germ of a curve, with the logarithmic structure pulled back from this germ. One then attempts to smooth this curve together with the map. The obstructions to deformations of a logarithmic map $[f:C \to X(\Delta)]$ are controlled by the group $H^1(C,f^\star T^{\mathrm{log}}X(\Delta))$,  and the logarithmic tangent bundle of $X(\Delta)$ is $\mathscr O_{X(\Delta)}^{\dim X}$. Thus, when $C$ has positive arithmetic genus, maps may be obstructed. However, after composing with the canonical map $X(\Delta)\to \cA(\Delta)$, we obtain a logarithmic map to a stack which is logarithmically \'etale, i.e. $T^{\mathrm{log}}\cA(\Delta) = 0$. Intuitively, since $H^1(\cA(\Delta),T^{\mathrm{log}}\cA(\Delta))$ is now trivial, one expects deformations of the map to be unobstructed. However, in order to make this precise, we would have to rely on results that have not appeared in the literature. Namely, the deformation theory of non-representable logarithmically \'etale morphisms. Instead, we appeal to Abramovich and Wise's study in~\cite{AW} of the moduli space of logarithmic prestable maps to Artin fans. This approach has the added benefit of allowing tropical curves whose edge lengths that are not rational.

\subsection{The stack of minimal logarithmic maps to the Artin fan} Logarithmic (pre)-stable maps were introduced in the papers~\cite{AC11,Che10,GS13}. In~\cite{AW}, an algebraic stack of \textit{minimal} logarithmic pre-stable maps to $\cA(\Delta)$ is constructed, and denoted $\fM(\cA(\Delta))$. The notion of minimality can be understood as follows. One wishes to work with the object $\fM(\cA(\Delta))$ as a moduli stack over the category of schemes, rather than over logarithmic schemes. Thus, a map from a test scheme $\underline S\to \fM(\cA(\Delta))$ should parametrize families of logarithmic pre-stable maps over $\underline S$. However, in order to build such a family, it is necessary to give the base scheme $\underline S$ a logarithmic structure. A priori there are numerous logarithmic structures that one may place on $\underline S$. A major insight in~\cite{Che10,GS13} is that there are distinguished \textit{minimal} logarithmic structures, which can be understood as the \textit{minimal requirements} that a logarithmic (pre)-stable map needs to satisfy. In other words, every other logarithmic stable map can be obtained from the minimal one by pulling back. Giving $\underline S$ this \textit{minimal} logarithmic structure, $\fM(\cA(\Delta))$ can be understood as a moduli stack parametrizing minimal logarithmic maps. We refer the reader to loc. cit. for further details.

\subsection{Step I: Curve and target} Fix a weighted balanced polyhedral complex $\mathscr P$ of dimension $1$ as before, and let $\Delta$ be any complete fan whose one skeleton supports the recession fan of $\mathscr P$. We may replace $\mathscr P$ with a subdivision such that each edge of $\mathscr P$ is completely contained in a single cone of $\Delta$. Consider a vertex $u$ of $\mathscr P$ incident to edges $e_1^u,\ldots, e_r^u$. Associate to $u$ a marked curve $C_u\cong \PP^1$, marked at distinct points $p_1^u,\ldots, p_r^u$, with the marked points in bijection with the edges emanating from $u$. Identify marked points $p_i^u$ and $p_j^w$ when the edges $e_i^u$ and $e_j^w$ coincide in $\mathscr P$, to form a nodal curve $\underline C_0$ whose dual graph coincides with $\mathscr P$. 

\subsection{Step II: A logarithmic map to $X(\Delta)$} Note that the linear constraints on the lengths of edges of a fixed combinatorial type are defined over $\QQ$, and thus, if one has a solution to these constraints, we can find edge lengths that are integral, to solve the same constraints. As a result, we may choose a polyhedral complex $\hat{\mathscr P}$ of the same combinatorial type as $\mathscr P$, but such that for every bounded edge $e$, the length $\ell(e)$ is an integral multiple of its weight $\omega(e)$. Let $u$ be a vertex of $\hat{\mathscr P}$, lying in the relative interior of a cone $\tau\in \Delta$. We now construct a logarithmic map $C_u\to V(\tau)$. The star around $u$ in $\mathscr P$ maps naturally to the fan of the toric variety $V(\tau)$. Perform a subdivision $\tau'\to \tau$ such that the star of $u$ maps onto the one skeleton of $\tau'$. Choose a map $f_u: C_u\to V(\tau')$ such that for any marked point $p$ corresponding to an edge $e$ of $\mathrm{star}(u)$, the point $p$ is mapped to the divisor $D_e$ of $V(\tau')$. Moreover, we require that the contact order of $f_u$ at $p$ is equal to the weight of $\mathscr P$ along $e$. If an edge of $\mathrm{star}(u)$ is contracted in $\tau'$, take the contact order to be $0$. Such maps can be explicitly written in the homogeneous coordinates of the toric variety\footnote{It is at this point that balancedness is being used.}, see for instance~\cite[Proposition 3.3.3]{R15b} for an explicit calculation. The map $f_u$ thus naturally has the structure of a logarithmic map. Composing with $V(\tau')\to V(\tau)$ and pushing forward the logarithmic structure as in~\cite[Appendix B]{AMW12}, we obtain a logarithmic map to $V(\tau)$. Ranging over all vertices $u$ and gluing, we obtain a map $f: \underline C_0\to X(\Delta)$. The curve $\underline C_0$ naturally acquires the structure of a logarithmically smooth curve over $\spec(\NN\to \CC)$, with the stalk of the characteristic of a node associated to an edge $e$ equal to the pushout $Q_e$
\[
\begin{tikzcd}
\NN \arrow{r}{\cdot \ell(e)/\omega(e)} \arrow[swap]{d}{\textrm{diag}} & \NN \arrow{d} \\
\NN^2 \arrow{r} & Q_e.
\end{tikzcd}
\]

We now promote $f$ to a logarithmic map over $\spec(\NN\to \CC)$. At the marked points, this data is equivalent to the contact orders described above. At the generic point of a component $C_u$ of $\underline C_0$, the stalk of the characteristic is $\NN$. Assume $C_u$ maps to the stratum corresponding to a cone $\delta\in \Delta$. To give a logarithmic map to $X(\Delta)$ here is equivalent to choosing a homomorphism $S_\delta\to \NN$. We choose this homomorphism to be the one given by the lattice point $u\in N$. Finally, consider a node $q$ corresponding to an edge $e$. The monoid $Q_e$ above can alternatively be described as~\cite[Remark 1.2]{GS13},
\[
Q_e = \{(n_1,n_2)\in \NN\times \NN: n_2-n_1 \in \ell(e)/\omega(e)\ZZ\}.
\]
Let $M_q$ be the dual lattice of the stratum to which $q$ maps. Following~\cite[Construction 1.16]{GS13}, to give a map $f^\flat_q:M_q\to Q_e$ is equivalent to the data of a homomorphism $u_q: M_q\to \ZZ$, such that 
\[
(n_2-n_1)\circ f^\flat_q(m) = u_q(m)\cdot \rho_q,
\]
together with the data of vertices associated to the components meeting at $q$. Equivalently, we may dualize, and set $u_q$ equal to the quantity $(v^q_1-v^q_2)/\rho_q$, where $v^q_1$ and $v^q_2$ are the vertices corresponding to the components meeting at $q$ and $\rho_q$ is the weight on the edge corresponding to $q$. We have thus produced a logarithmic map $f:C_0\to X(\Delta)$ over $\spec(\NN\to \CC)$.

\subsection{Step {III}: Working over the minimal base} By composing the map constructed in the previous step with the quotient $X(\Delta)\to \cA(\Delta)$, we obtain a map $C_0\to \cA(\Delta)$. By definition of the moduli space $\fM(\cA(\Delta))$, this is equivalent to a moduli map
\[
\spec(\NN\to \CC) \to \spec(Q\to \CC)\to \fM(\cA(\Delta)),
\]
where $Q$ is the stalk of the minimal characteristic of $\fM(\cA(\Delta))$ at the image of the underlying map $\spec(\CC)\to {\fM(\cA(\Delta))}$ given by the composition. We now discard the map on the left and consider only the second arrow. It is proved in~\cite[Proposition 1.5.1]{AW} that the stack of maps ${\fM(\cA(\Delta))}$ is logarithmically smooth, and as a consequence, we obtain a family of maps as below:
\begin{equation}\label{eqn: family-of-maps}
\begin{tikzcd}
\widetilde{\mathscr C} \arrow{r} \arrow{d} & \mathscr C^{\textrm{univ}} \arrow{r}\arrow{d} & \cA(\Delta) \\
\spec(Q\to \CC\llbracket Q\rrbracket) \arrow{r} & \fM(\cA(\Delta)). &
\end{tikzcd}
\end{equation}

\subsection{Step IV: Constructing the analytic map} Our next step is to observe an explicit description of the minimal characteristic $Q$. The characteristic used in~\cite{AW} follows the definition of Abramovich and Chen~\cite{AC11}. However, by~\cite[Proposition 4.8]{AC11}, this minimal characteristic coincides with the basic monoid of Gross and Siebert. This follows from the universal properties of the basicness and minimality conditions. In turn, by~\cite[Remark 1.21]{GS13}, the minimal characteristic $Q$ above coincides with the dual monoid of the monoid of integral points in the cone of tropical curves in $\Delta$ having combinatorial type the same as $\mathscr P$. 

Summarizing, the cone of tropical curves having combinatorial type $\mathscr P$ is precisely $\Hom(Q,\RR_{\geq 0})$, and thus $\mathscr P$ defines a monoid homomorphism $Q\to \RR_{\geq 0}$. Combining this with the final diagram in the above step, we have a morphism
\[
\spec(\CC\llbracket \RR_{\geq 0}\rrbracket) \to \spec(\CC\llbracket Q \rrbracket) \to \fM(\cA(\Delta)).
\]
By pulling back the universal curve and map, we have
\[
\begin{tikzcd}
\mathscr C \arrow{r}\arrow{d} & \cA(\Delta) \\
\spec(\CC\llbracket \RR_{\geq 0}\rrbracket). & \\ 
\end{tikzcd}
\]
View the ring $\CC\llbracket \RR_{\geq 0}\rrbracket$ as a generalized power series ring with the standard valuation, and let $C$ be the generic fiber of $\mathscr C$. Apply Raynaud's generic fiber construction to obtain a map of analytic stacks $\varphi: C^{\an}\to [X^{\an}/T^{\an}_\circ]$. 

\subsection{Step V: Computing the tropicalization} It remains to show that the image of this map $\varphi$ is the polyhedral complex $\mathscr P$. To do this, return to the family of maps over the base $\spec(Q\to \CC\llbracket Q\rrbracket)$ constructed in Step III (\ref{eqn: family-of-maps}). Consider this base as a scheme over the trivially valued field $\CC$, and apply the functor $(\cdot)^\beth$. Note that $Q$ is a toric monoid, and hence $\spec(\CC\llbracket Q\rrbracket)$ is the local model appearing in Thuillier's construction of extended skeletons for formal fibers of toroidal embeddings over the trivial valuation. The formation of extended skeletons for $\beth$-spaces is functorial~\cite[Theorem 1.1]{U13} so we obtain
\[
\begin{tikzcd}
(\widetilde{\mathscr C})^\beth \arrow{r} \arrow{d} & {|\cA(\Delta)\times \spec(\CC\llbracket Q\rrbracket)^\beth|} \arrow{d} \\
\overline \Sigma(\widetilde{\mathscr C}) \arrow{r}  & N(\Delta) \times \Hom(Q,\overline \RR_{\geq 0}).
\end{tikzcd}
\]
Here, $\overline \Sigma(-)$ is the compactified cone complex associated to a toroidal embedding. We briefly recall that if $\sigma = \Hom(S_\sigma, \RR_{\geq 0})$ is a cone, then its compactification is given by $\overline \sigma = \Hom(S_\sigma, \RR_{\geq 0}\sqcup \{\infty\})$. The compactified cone complex is obtained by gluing the compactified cones associated to the toroidal embedding. See~\cite{Thu07,U13} for details concerning its precise definition. 

The tropical curve gives rise to a point $[\mathscr P]\in \Hom(Q,\RR_{\geq 0})$, and moreover there is a natural map 
\[
\bm p: \overline \Sigma(\mathscr C')\to \Hom(Q,\overline \RR_{\geq 0}).
\]
The image of $\varphi$ coincides with the image of the slice $\bm p^{-1}([\mathscr P])\subset \overline \Sigma(\mathscr C')$ in $N(\Delta)\times \{[\mathscr P]\}$. By construction, the image of $\bm p^{-1}([\mathscr P])\subset \overline \Sigma(\mathscr C')$ has the correct combinatorial type, so we need only check the lengths of its bounded edges. But observe that the length of an edge $e$ is equal to the valuation of the deformation parameter $f_q$ of the corresponding node $q$. We have chosen the valuation in the previous step to be the one that assigns to $f_q$ valuation $\ell(e)/\omega(e)$. The result follows. \qed

\subsection*{Acknowledgements} I learned about the idea that tropical curves should be
realizable by maps to the Artin fan from Dan Abramovich~\cite{A14-Slides}, and I thank
him for encouraging me to pursue the problem further, as well as for comments on previous drafts. I am grateful to Dori
Bejleri, Tyler Foster, Sam Payne, David Speyer, Martin Ulirsch, and Jonathan Wise for
illuminating conversations, and to David Speyer for his permission to use Figure~\ref{speyercurve}. I am especially grateful to both referees for their careful reading, comments, and corrections. This work was initiated when I was a visiting student at Brown University in Spring 2015, and I am grateful to the Department of Mathematics for ideal working conditions.

\bibliographystyle{siam}
\bibliography{Superabundance}

\end{document}